\documentclass{aupl}

\usepackage{
amsfonts,
latexsym,
amssymb,
amsmath,
amsthm,
enumerate,
verbatim,
mathrsfs,
epigraph,
}

\usepackage{url}

\newcommand{\labbel}[1]{\label{#1} [[{\bf #1}]]}  
\renewcommand{\labbel}{\label}

\newtheorem{theorem}{Theorem}[section]
\newtheorem{lemma}[theorem]{Lemma}

\theoremstyle{definition}
\newtheorem{definition}[theorem]{Definition}

\newtheorem{problem}[theorem]{Problem} 
\newtheorem{problems}[theorem]{Problems} 

\theoremstyle{remark}

\numberwithin{equation}{section}

 \allowdisplaybreaks[4]

\begin{document}

\title{Best values for the distributivity spectrum}

\author{Paolo Lipparini} 
\address{Dipartimento Quaternario di Matematica\\Viale della  Ricerca
 Scientifica\\Universit\`a di Roma ``Tor Vergata'' 
\\I-00133 ROME ITALY}
\urladdr{http://www.mat.uniroma2.it/\textasciitilde lipparin}


\thanks{Work performed under the auspices of G.N.S.A.G.A. 
The author acknowledges the MIUR Department Project awarded to the
Department of Mathematics, University of Rome Tor Vergata, CUP
E83C18000100006.}

\begin{abstract}
We provide optimal bounds for
$\alpha (\beta \circ \gamma \circ  \dots )$
in $4$-distributive varieties, as well as some further partial
generalizations. 
\end{abstract} 

\maketitle

\section{Introduction} \labbel{int} 

B. J{\'o}nsson \cite{JD} showed that
a variety $\mathcal V$ is congruence distributive if and only if
there is some $n$ and there are terms
$t_0, \dots, t_n$ satisfying an appropriate set of
equations. If this holds for some
specific $n$, then $\mathcal V$ is said to be
$n$-distributive. Varieties 
which are $n$-distributive for  
small values of $n$ satisfy  special properties.
The case $n=2$ means the existence
of a majority term.
Special results holding for $n=3$
and $n=4$ are presented in 
\cite{baker,CDMM,KV,jds,ia,B,daysh}.
See also \cite{Fi,FV,adjt,KaV,Ke,Le} 
 for more general results about ``distributivity levels''.
 
In \cite{jds} we showed that if $\mathcal V$ is 
a $k{+}1$-distributive variety, then, for every $m \geq 1$, $\mathcal V$   
satisfies the congruence identity
\begin{equation}\label{m}  
 \alpha ( \beta \circ \gamma \circ   {\stackrel{m+1}{\dots}}  ) \subseteq 
\alpha \beta \circ \alpha \gamma \circ  {\stackrel{km+1}{\dots}}    
\end{equation}      
In identities like \eqref{m} juxtaposition denotes intersection and
$\beta \circ \gamma \circ {\stackrel{\ell}{\dots}} $
is $\beta \circ \gamma \circ \beta \circ \gamma \circ \dots $
with $\ell-1$ occurrences of  $ \circ $.  
When we know that, say, $\ell$ is even, we sometimes write
$\beta \circ \gamma \circ {\stackrel{\ell}{\dots}} \circ \gamma  $
when we want to make clear that the last factor is $\gamma$.
Other unexplained notions and notations are from
\cite{jds,daysh}.

In particular, it follows from \cite{jds}
that every $4$-distributive variety satisfies  
\begin{equation}\label{m4}  
 \alpha ( \beta \circ \gamma \circ   {\stackrel{m+1}{\dots}}  ) \subseteq 
\alpha \beta \circ \alpha \gamma \circ  {\stackrel{3m+1}{\dots}}   
\end{equation}      
for every $m \geq 1$.
We show that in the special case of $4$-distributive varieties
\eqref{m} and \eqref{m4}  can be improved.

\begin{theorem} \labbel{gen}
Every $4$-distributive variety
satisfies the congruence identities 
\begin{align} \labbel{1} 
\alpha (\beta \circ \gamma \circ  {\stackrel{n}{\dots}} \circ \gamma)
&\subseteq 
  \alpha \beta \circ \alpha \gamma
 \circ {\stackrel{2n}{\dots}} \circ     \alpha \gamma,
&&\text{for $n$ even,}
\\
\labbel{2}
\alpha (\beta \circ \gamma \circ  {\stackrel{n}{\dots}} \circ \beta )
&\subseteq 
  \alpha \beta \circ \alpha \gamma
 \circ {\stackrel{2n+1}{\dots}} \circ     \alpha \beta,
&&\text{for $n$ odd.} 
 \end{align}    
 \end{theorem} 

Theorem \ref{gen} shall be proved at the end of Section \ref{e4}.  
Results hinted in 
 \cite[Remark 10.1(a)(1)-(3)]{daysh} imply that the bounds are optimal.

\section{A simple introductory example} \labbel{simpl} 

We first present a special case of our main result.
The proof of the next theorem is relatively simple, but it
contains most of the ideas we shall use in the general case.
However, many technical details will be needed to 
prove the general results in the next section, hence this section
might serve as a gentle introduction to our methods. We first
recall the definition of J{\'o}nsson terms in the special case $n=4$.

\begin{definition} \labbel{4def}
A variety $\mathcal V$ is \emph{$4$-distributive}  
or \emph{$4$-J{\'o}nsson} if
it has terms $t_1, t_2, t_3$
such that the equations
\begin{align*} \labbel{cd4}  
x&= t_i(x,y,x), & \text{ for  } i=1, 2, 3;
\\
 x&=t_1(x,x,z),   & t_1(x,z,z)&=t_2(x,z,z), 
\\
t_2(x,x,z) &=t_3(x,x,z),   &    t_3(x,z,z)&=z, 
 \end{align*} 
are valid in $\mathcal V$.

If we do not require the equation
$x=t_2(x,y,x)$, namely, if the first line
is assumed to hold only for $i=1$ and  $i=3$,
we will speak of   a \emph{defective $4$-J{\'o}nsson variety}.
 \end{definition}    

A variety $\mathcal V$ is $4$-distributive if and only if $\mathcal V$  
satisfies the congruence identity $ \alpha ( \beta \circ \gamma ) 
\subseteq \alpha \beta \circ \alpha \gamma \circ \alpha \beta \circ \alpha \gamma$.
See \cite[Section 2]{daysh} for more information.  
Similarly, a variety is defective $4$-J{\'o}nsson
if and only if it satisfies
the congruence identity $ \alpha ( \beta \circ \gamma ) 
\subseteq \alpha \beta \circ \alpha (\gamma \circ \beta) \circ \alpha \gamma$.

Notice that every congruence permutable variety is 
defective $4$-J{\'o}nsson: just take $t_1$ and  $t_3$
to be projections. On the other hand, Polin's variety 
is  defective $4$-J{\'o}nsson but not congruence
 modular \cite[Remark 10.11]{daysh}.

Frequently, Definition \ref{4def} is given by considering
further terms $t_0$ and  $t_4$ which are the trivial projections
onto the first, respectively, the third coordinate.   
Then the condition $x=t_1(x,x,z)$ 
is written as $t_0(x,x,z)=t_1(x,x,z)$ and symmetrically at the end.
This alternative convention justifies the name ``$4$-distributive''.
However, here it will be simpler to deal with $t_1, t_2, t_3$ alone.

\begin{theorem} \labbel{4d}
Every $4$-distributive variety
satisfies the congruence identity 
\begin{equation*}
\alpha (\beta \circ \gamma \circ \beta \circ \gamma )
\subseteq 
  \alpha \beta \circ \alpha \gamma
 \circ {\stackrel{8}{\dots}} \circ     \alpha \gamma.
 \end{equation*}    
 \end{theorem} 

\begin{proof}  
Suppose that $\mathcal V$ is $4$-distributive with  J{\'o}nsson terms
$t_1$, $t_2$, $t_3$.
If $(a,e) \in \alpha (\beta \circ \gamma \circ \beta \circ \gamma )$,
then $ a \mathrel { \alpha } e $ and $ a \mathrel { \beta } b \mathrel { \gamma } c 
\mathrel { \beta } d \mathrel { \gamma } e  $,
for some elements $b,c,d$ of the algebra under consideration.
Let $F= t_2(t_2(a,d,e),t_2(a,c,e),t_2(a,b,e))$.  Then
\begin{align*}
a = t_1(a,a,t_2(a,c,e)) =& t_1(t_1(a,a,e),t_2(a,a,a),t_2(a,c,e))  \mathrel \beta    
\\
& t_1(t_1(a,b,e),t_2(a,b,b),t_2(a,c,e))  \mathrel \gamma   
\\
& t_1(t_1(a,c,e),t_2(a,c,c),t_2(a,c,e))  \mathrel \beta    
\\
& t_1(t_1(a,d,e),t_2(a,c,d),t_2(a,c,e))  \mathrel \gamma   
\\
& t_1(t_1(a,e,e),t_2(a,c,e),t_2(a,c,e)) =
\\
& t_2(t_2(a,e,e),t_2(a,c,e),t_2(a,c,e)) \mathrel { \gamma } 
\\
& t_2(t_2(a,d,e),t_2(a,c,e),t_2(a,b,e)) = F.
\end{align*}    

Since $a \mathrel { \alpha } e $,
all the elements in the above list are $\alpha$-equivalent to $a$,
by the equations $t_i(x,y,x)= x$, hence 
$(a, F) \in \alpha \beta \circ \alpha \gamma \circ \alpha \beta \circ \alpha \gamma  $.
Symmetrically,
$(F,e) \in \alpha \beta \circ \alpha \gamma \circ \alpha \beta \circ \alpha \gamma  $,
thus
$(a, e) \in \alpha \beta \circ \alpha \gamma
 \circ {\stackrel{8}{\dots}} \circ     \alpha \gamma  $.
\end{proof}  

In the terminology from \cite{jds},
Theorem \ref{4d} asserts that $ J _{ \mathcal V} (3) \leq 7 $,
for every $4$-distributive variety $\mathcal V$.    
As we mentioned, this improves the bound $ J _{ \mathcal V} (3) \leq 9$
provided by \cite[Corollary 2.2]{jds} in the case 
$m=1$ and $k= \ell=3$.

Theorem \ref{4d} provides an optimal result, as shown
by \cite[Theorem 2.1]{B}, case $n=4$, equation (2-3).

The proof of Theorem \ref{4d}
applies to the case of 
defective $4$-J{\'o}nsson varieties,
noticing that the term $t_2$ is idempotent in any case,
but we need the further assumption $ a \mathrel { \alpha  } c$.

\begin{theorem} \labbel{4dd}
Every defective $4$-J{\'o}nsson variety
satisfies the congruence identity 
\begin{equation*}
 \alpha ( \beta \circ \gamma ) \circ \alpha ( \beta \circ \gamma ) 
\subseteq \alpha \beta \circ \alpha \gamma \circ
\alpha \beta \circ \alpha (\gamma \circ  
 \beta) \circ \alpha \gamma \circ
\alpha \beta \circ \alpha \gamma. 
 \end{equation*}    
 \end{theorem}

\section{Exact bounds for $4$-distributive varieties} \labbel{e4}

Throughout, we fix some algebra $\mathbf M$, 
congruences $\alpha, \beta, \gamma $ on $\mathbf M$,
$n , \ell \in \mathbb N$, $n, \ell $ even, $n>0$   and elements
\begin{equation*} 
a=b_0 \mathrel{ \beta } b_1 \mathrel{ \gamma  } b_2
 \mathrel{ \beta } b_3 \mathrel{ \gamma  } b_4 \dots 
b_{n-1}\mathrel{ \gamma  } b_n =c  
 \end{equation*}     
 of $M$
such that $a \mathrel { \alpha  } c$
(here $c$ plays the role played by $e$ in the previous section).
An \emph{$n$-$\ell$-$n$-chain}  is a sequence
$A_0 ,  A_1, \dots, \allowbreak  A_{n-1}, 
B_0, \dots, \allowbreak  B_\ell ,  C_1 , \dots,  C_n$ 
of elements of $M$ 
such that
\begin{multline*} 
a=A_0 \mathrel{ \beta } A_1
\mathrel{ \gamma  } A_2 \mathrel { \beta  } \dots
  \mathrel{ \gamma    } A_{n-2} \mathrel{ \beta   } A_{n-1}
\mathrel { \gamma  } B_0
\mathrel{ \beta  } B_1 \mathrel{ \gamma  } B_2 \dots
\\
\mathrel{ \beta  } B_{\ell-1}  \mathrel { \gamma  } B_\ell 
 \mathrel{ \beta } C_1
\mathrel{ \gamma  } C_2 \mathrel { \beta }  \dots \mathrel{ \gamma  }  C_n=c
 \end{multline*}
(to make the chain more symmetrical, we could consider
$A_n=B_0$ and $C_0= B_\ell$, but we shall not need this)
and, furthermore,
  \begin{enumerate}[(C1)]   
 \item 
All the elements $A_i$, $B_j$, $C_k$ are pairwise $ \alpha $-related, and    
\item
For every $j \leq \ell$,
both $(a,B_j)$ and $(B_j, c)$ belong to 
$ \beta  \circ  \gamma \circ  
 {\stackrel{n}{\dots}} \circ  \gamma    $. 
  \end{enumerate}  
In particular, for an $n$-$0$-$n$-chain, 
we have $B_0=B_\ell$,
thus in this case $ (a, c) \in \alpha \beta  \circ \alpha \gamma \circ  
 {\stackrel{2n}{\dots}} \circ \alpha \gamma    $.

\begin{lemma} \labbel{lemch}
If $\mathbf M$ belongs to a $4$-distributive variety,
then, for all  $\alpha, \beta, \gamma $, even $n$
and $a, b_1, \dots $ as above, 
 there is an $n $-$0$-$n$-chain. 
 \end{lemma} 

\begin{proof}
The case $n=2$ is J{\'o}nsson's argument.  
In general, the same argument provides an
$n $-$(n-2)$-$n$-chain. Indeed, this is witnessed by the elements
\begin{align*}
a=A_0 & =t_1(a,a,c) 
\\
A_1 & =t_1(a,b_1,c) 
\\
&\dots
\\
A_i & =t_1(a,b_i,c)
\\
&\dots
\\
A_{n-1} & =t_1(a,b_{n-1},c)
\\
B_0 & =t_2(a,b_{n-1},c) 
\\
B_1 & =t_2(a,b_{n-2},c) 
\\
&\dots
\\
B_{n-2} & =t_2(a,b_{1},c)
\\
C_1 & =t_3(a,b_{1},c) 
\\
&\dots
\\
C_i & =t_3(a,b_i,c)
\\
&\dots
\\
C_{n-1} & =t_3(a,b_{n-1},c)
\\
C_{n} & =t_3(a,c,c) =c
  \end{align*} 
since 
$A_{n-1} =t_1(a,b_{n-1},c) \mathrel { \gamma  } 
t_1(a,c,c)= t_2(a,c,c) \mathrel { \gamma  } t_2(a,b_{n-1},c) = B_0 $
and symmetrically  $B_{n-2} \mathrel { \beta  } C_1$. 
The only nontrivial thing which remains to be checked is (C2).
Indeed, for every $j$, say, $j$ odd, 
$a=t_2(a,a,a) \mathrel { \beta  } t_2(a,b_1,b_1)
 \mathrel { \gamma   } t_2(a,b_2,b_2) \mathrel { \beta  }
\dots  \mathrel { \gamma  }   t_2(a,b_{n-1-j},b_{n-1-j}) \mathrel { \beta   } 
t_2(a,b_{n-1-j},b_{n-j}) \mathrel { \gamma  }
 t_2(a,b_{n-1-j},b_{n-j+1})  \dots t_2(a,b_{n-1-j},b_{n-1})
\mathrel { \gamma  } t_2(a,b_{n-1-j},c) = B_j $.  
Notice that we are not requiring that the elements witnessing
(C2) are $\alpha$-related.

In order to complete the proof it is then enough to show   
that if $\ell >0$, $\ell$ is even and there is an 
 $n $-$\ell$-$n$-chain, then there is an 
$n $-$(\ell - 2)$-$n$-chain. The conclusion follows
from a finite induction.

So let $A_0 ,  \dots, B_0, \dots ,   C_n$ 
be an $n $-$\ell$-$n$-chain.
By (C2), there are $X_i$s such that
$a=X_0 \mathrel { \beta  } X_1 \mathrel { \gamma  }
\dots \mathrel { \beta  } X_{n-1} \mathrel { \gamma  }
X_n=  B_{\ell-1}$ and, symmetrically,
$Y_i$s such that
$B_1=Y_0 \mathrel { \beta  } Y_1 \mathrel { \gamma  }
\dots \mathrel { \beta  } Y_{n-1} \mathrel { \gamma  }
X_n= c$. 
Define
\begin{align*}
a=A'_0 & =t_1(a,a,B_{\ell-1}) 
\\
A'_1 & =t_1(A_1,X_1,B_{\ell-1}) 
\\
&\dots
\\
A'_i & =t_1(A_i,X_i,B_{\ell-1})
\\
&\dots
\\
A'_{n-1} & =t_1(A_{n-1},X_{n-1},B_{\ell-1})
\\
B'_0 & =t_2(B_0,B_{\ell-1},B_{\ell})
\\
B'_1 & =t_2(B_0,B_{\ell-2},B_{\ell})
\\
&\dots
\\
B'_{\ell-2} & =t_2(B_0,B_{1},B_{\ell})
\\
C'_1 & =t_3(B_1,Y_{1},C_1) 
\\
&\dots
\\
C'_i & =t_3(B_1,Y_{i},C_i)
\\
&\dots
\\
C'_{n-1} & =t_3(B_1,Y_{n-1},C_{n-1})
\\
C'_{n} & =t_3(B_1,c,c) =c
  \end{align*} 
This is an $n $-$(\ell - 2)$-$n$-chain. Indeed, 
\begin{align*}  
A'_{n-1}&=t_1(A_{n-1},X_{n-1},B_{\ell-1}) \mathrel { \gamma  }
 t_1(A_{n-1},B_{\ell-1},B_{\ell-1})
\\
& = t_2(A_{n-1},B_{\ell-1},B_{\ell-1}) 
\mathrel { \gamma  } t_2(B_0,B_{\ell-1},B_{\ell}) = B'_0
\end{align*} 
and, symmetrically, $B'_{\ell-2} \mathrel { \beta  } C'_1$.

Notice that the case $\ell=2$ is allowed in the above computations; in such a
case we have  $B'_0=B'_{\ell-2}=t_2(B_0,B_{1},B_{2})$ and we get
$A'_{n-1}
\mathrel { \gamma  } B'_0
 \mathrel{ \beta } C'_1$  
thus we actually reach an $n $-$0$-$n$-chain
at the end of the induction. Compare the different
situation when $n$ is odd, which
shall be treated below.  

It remains to check that the new chain 
$A'_0 ,  \dots, B'_0, \dots ,   C'_n$
satisfies (C2). This follows from
$a=t_2(a,a,a)$ and the assumption about the
original chain that   
$(a,B_0)$, $(a,B_{\ell-j-1})$ and $(a,B_\ell)$ all belong to
$ \beta  \circ  \gamma \circ  
 {\stackrel{n}{\dots}} \circ  \gamma    $. 
Indeed, suppose that the above relations are witnessed by sequences
$a=f_0 \mathrel { \beta  } f_1 \dots $,
$a=g_0 \mathrel { \beta  } g_1 \dots $ and
$a=h_0 \mathrel { \beta  } h_1 \dots $,
then 
$a=t_2(a,a,a) \mathrel { \beta  } t_2(f_1,g_1,h_1) 
\mathrel { \gamma  } t_2(f_2,g_2,h_2) \dots
t_2(B_0,B_{\ell-j-1},B_\ell)= B_j$
witnesses (C2) for the new chain 
$A'_0 ,  \dots, B'_0, \dots ,   C'_n$.
\end{proof}

   \begin{proof}[Proof of Theorem \ref{gen}]
The inclusion \eqref{1} is immediate from Lemma \ref{lemch}
and the comment immediately preceding it.

In order to prove \eqref{2}, we should modify the definition of an     
$n$-$\ell$-$n$-chain in the case when $n$ and $\ell$ are odd.
The definition is very similar, but in this case we have
  $\dots \mathrel{ \gamma } A_{n-1}
\mathrel { \beta  } B_0
\mathrel{ \gamma } B_1 \mathrel{ \beta  } B_2 \dots
 \mathrel { \gamma  }  B_{\ell-2} \mathrel{ \beta  }
 B_{\ell-1}  \mathrel { \gamma  } B_\ell 
 \mathrel{ \beta } C_1 \mathrel { \gamma  } \dots $.
Since here $n$ is odd, condition (C2) now reads:
  \begin{enumerate} 
   \item[(C2, $n$ odd)]  
both $(a,B_j)$ and $(B_j, c)$ belong to 
$ \beta  \circ  \gamma \circ  
 {\stackrel{n}{\dots}} \circ  \beta   $. 
  \end{enumerate}
 
In the reduction procedure in the analogue of the proof
of Lemma \ref{lemch} we should consider instead
\begin{align*}
&\dots
\\
A'_{n-2} & =t_1(A_{n-2},X_{n-2},B_{\ell-1})
\\
A'_{n-1} & =t_1(A_{n-1},X_{n-1},B_{\ell-1})
\\
B'_0 & =t_2(B_0,B_{\ell-2},B_{\ell-1})
\\
B'_1 & =t_2(B_0,B_{\ell-3},B_{\ell})
\\
B'_2 & =t_2(B_0,B_{\ell-4},B_{\ell})
\\
&\dots
\\
B'_{\ell-2} & =t_2(B_0,B_{0},B_{\ell})=t_3(B_0,B_{0},B_{\ell})
\\
C'_1 & =t_3(B_0,Y_{1},C_1) 
\\
C'_2 & =t_3(B_0,Y_{2},C_2) 
\\
&\dots,
  \end{align*} 
which is an $n$-$(\ell-2)$-$n$-chain, since
\begin{align*}  
A'_{n-1} & = t_1(A_{n-1},X_{n-1},B_{\ell-1}) \mathrel { \beta  }
 t_1(A_{n-1},B_{\ell-1},B_{\ell-1})
\\
 &= t_2(A_{n-1},B_{\ell-1},B_{\ell-1})
\mathrel { \beta  } t_2(B_0,B_{\ell-2},B_{\ell-1})=B'_0.
\end{align*} 

In comparison with the case $n$ and $\ell$ even, notice that
here we cannot proceed with the induction when $\ell=1$,
since in the definition of $B'_0$ we have used $B_{\ell-2}$,
which is undefined when $\ell=1$. This explains the shift by $1$
in \eqref{2} with respect to \eqref{1}.       
 \end{proof}

\section{Further remarks} \labbel{fur}

We are not claiming that the next problems are difficult.

\begin{problems} \labbel{probs}
(a) Can we improve the conclusion of Theorem \ref{4dd} to
\begin{equation*}
 \alpha ( \beta \circ \gamma \circ  \beta \circ \gamma ) 
\subseteq \alpha \beta \circ \alpha \gamma \circ
\alpha \beta \circ \alpha (\gamma \circ  
 \beta) \circ \alpha \gamma \circ
\alpha \beta \circ \alpha \gamma \ ? 
 \end{equation*}    

(b) Devise a generalization of Theorem \ref{4dd}
along the lines of (Lemma \ref{lemch} and) Theorem \ref{gen}.  

(c) Generalize the present results when
$\alpha$, $\beta$ and $\gamma$ are only assumed to 
be reflexive and admissible relations. Compare some results
in \cite{jds,B}. More generally, consider
the other spectra introduced in \cite[Section 3]{jds}.  
Some more spectra, possibly for non distributive varieties,
are implicit in \cite{daysh}.
 \end{problems}

We have seen that, in general, the results from
\cite{jds} are not the best possible ones.
However, the present discussion deals only with
$4$-distributive varieties (for which we have indeed found optimal bounds).

\begin{problem} \labbel{prob} \cite{jds} 
For $n>4$ and $m>1$, find the smallest possible value of $k$
such that every $n$-distributive variety satisfies
\begin{equation}\label{mk}  
 \alpha ( \beta \circ \gamma \circ   {\stackrel{m+1}{\dots}}  ) \subseteq 
\alpha \beta \circ \alpha \gamma \circ  {\stackrel{k}{\dots}}    
\end{equation}
 \end{problem}

We now present some partial improvements on Theorem \ref{4d}.
For $m=3$ and $k=n-1$, equation \eqref{m} reads
\begin{equation} \labbel{mu}    
 \alpha ( \beta \circ \gamma \circ \beta \circ \gamma  ) \subseteq 
\alpha \beta \circ \alpha \gamma \circ  {\stackrel{3n-2}{\dots}}   
\end{equation}      
We now generalize Theorem \ref{4d} by showing that \eqref{mu}
can be improved by $2$.  We do not know whether the argument 
can be iterated in order to improve the bound.

\begin{theorem} \labbel{nd}
If $n\geq 4$ and $n$ is even, then every $n$-distributive variety
satisfies the congruence identity 
\begin{equation*}
\alpha (\beta \circ \gamma \circ \beta \circ \gamma )
\subseteq 
 \alpha \beta \circ \alpha \gamma
 \circ {\stackrel{3n-4}{\dots}} \circ     \alpha \gamma.
 \end{equation*}    
 \end{theorem} 

\begin{proof}  
In the present case, besides $t_1, t_2, t_3 $,
we have further terms $t_3, t_4, \dots$ satisfying 
similar identities.
Extend the chain from the proof of Theorem \ref{4d}
as follows.   To save some lines, we write, say,
$t_{2=3}(c,c,e)$ for an expression which can be,
equivalently, $t_{2}(c,c,e)$ or  $t_{3}(c,c,e)$, by an equation 
in Definition \ref{4def}.   
\begin{align*}
F&= t_2(t_2(a,d,e),t_2(a,c,e),t_2(a,b,e))   \mathrel { \beta}
\\
& t_2(t_2(a,c,e),t_2(a,c,e),t_2(a,a,e)) =
\\
& t_3(t_2(a,c,e),t_2(a,c,e),t_3(a,a,e))  \mathrel { \beta}
\\
& t_3(t_2(a,c,e),t_2(b,c,e),t_3(a,b,e))  \mathrel { \gamma } 
\\
& t_3(t_2(a,b,e),t_{2=3}(c,c,e),t_3(a,c,e))  \mathrel { \beta  }
\\
& t_3(t_{2=3}(a,a,e),t_3(c,d,e),t_3(a,d,e))  \mathrel { \beta } 
\\
& t_3(t_3(a,b,e),t_3(c,d,e),t_3(a,d,e))  \mathrel { \gamma  } 
\\
& t_3(t_3(a,c,e),t_{3=4}(b,e,e),t_{3=4}(a,e,e))   \mathrel { \gamma  }
\\
& t_3(t_3(a,c,e),t_4(b,d,e),t_4(a,d,e))  \mathrel { \beta }
\\
& t_{3=4}(t_3(a,d,e),t_4(a,c,e),t_4(a,c,e))  \mathrel { \gamma  }
\\
& t_4(t_{3=4}(a,e,e),t_4(a,c,e),t_4(a,c,e)) \mathrel { \gamma  }
\\
& t_4(t_4(a,d,e),t_4(a,c,e),t_4(a,b,e))
\\ 
 &\dots 
\end{align*}    
Iterate the above construction,
by steps of $2$, and at the end proceed in a symmetrical way
with respect to the proof of Theorem \ref{4d}. 
\end{proof}  

The next theorem slightly improves the case $\ell=4$
(or possible bilateral variations) of 
the alvin case 
in \cite[Corollary 6(2)(3)]{gumalv}.  

\begin{theorem} \labbel{ndalv}
If $n\geq 4$ and $n$ is even, then every $n+2$-alvin variety
satisfies the congruence identity 
\begin{equation*}
\alpha (\beta \circ \gamma \circ \beta \circ \gamma )
\subseteq 
 \alpha \gamma  \circ \alpha \beta 
 \circ {\stackrel{3n-2}{\dots}} \circ     \alpha \beta .
 \end{equation*}    
 \end{theorem} 

\begin{proof}
As an unusual convention, let us denote by $t_0$
the first term in the alvin sequence, thus $x=t_0(x,y,y)=t_0(x,y,x)$
and $t_0(x,x,z)=t_1(x,x,z)$ are valid equations.
This unusual convention allows us to keep the same indices as 
in the preceding parts of the manuscript.
By the proof of
\cite[Corollary 6(1)]{gumalv}, the alvin terms can be chosen in such a way that  
 $a \mathrel { \alpha \gamma {\circ} \alpha \beta   } t_0(a,a,e)$,
under the usual
assumption $(a,e) \in \alpha (\beta \circ \gamma \circ \beta \circ \gamma )$.
Thus
 $a \mathrel { \alpha \gamma } X \mathrel { \alpha \beta   } t_0(a,a,e)$,
for some $X$. Then
\begin{align*}
a= & t_0(a,t_1(a,b,e),t_1(a,b,e))  \mathrel {  \gamma   } 
\\
& t_0(X,t_1(a,b,e),t_1(a,c,e)) \mathrel { \beta  } 
\\
& t_0(t_0(a,a,e),t_1(a,a,e),t_1(a,c,e))    =
\\
& t_0(t_1(a,a,e),t_1(a,a,e),t_1(a,c,e)) = 
\\
& t_1(t_1(a,a,e),t_1(a,a,e),t_1(a,c,e))  \mathrel \beta    
\\
& t_1(t_1(a,b,e),t_1(a,b,e), t_1(a,d,e))  \mathrel \gamma   
\\
& t_1(t_1(a,b,e),t_1(a,b,e),t_1(a,e,e) = 
\\
& t_1(t_1(a,b,e),t_1(a,b,e),t_2(a,e,e)  \mathrel \gamma   
\\
& t_1(t_1(a,c,e),t_1(a,c,d),t_2(a,d,e))  \mathrel \beta   
\\
& t_1(t_1(a,c,e),t_1(a,c,c),t_2(a,c,e))  \mathrel =   
\\
& t_1(t_1(a,c,e),t_2(a,c,c),t_2(a,c,e))  \mathrel \beta    
\end{align*}    
and the last line in the above formula is equal to the 
third line in the displayed formula in the proof of Theorem \ref{4d},
hence we can apply the proofs of  Theorems \ref{4d} and \ref{nd}. 
\end{proof}

\end{document}